\documentclass[12pt,reqno]{amsart}
\usepackage{amssymb}
\usepackage{amsxtra}
\usepackage{amsmath}
\usepackage{amsfonts}

\numberwithin{equation}{section}

\newtheorem{thm}{Theorem}[section]

\newtheorem{prop}[thm]{Proposition}
\newtheorem{cor}[thm]{Corollary}

\newtheorem{lem}[thm]{Lemma}
\newtheorem{de}[thm]{Definition}
\newtheorem{rem}[thm]{Remark}
\newtheorem{ex}[thm]{Example}

\def\ty{{\tilde y}}
\def\w{\omega}
\def\bx{\bold x}
\def\av{\alpha^{\vee}}
\def\tw{\widetilde W}
\def\tv{\widetilde V}
\newcommand{\eqa}{\begin{eqnarray}}
\newcommand{\eeqa}{\end{eqnarray}}
\newcommand{\beq}{\begin{equation}}
\newcommand{\eeq}{\end{equation}}
\newcommand{\nn}{\nonumber}
\newcommand{\p}{\partial}

\newcommand{\pal}{\partial}

\newcommand{\td}{\tilde d}

\newcommand{\al}{\alpha}

\newcommand{\ta}{\theta}

\def \om{\mu}

 \def \dsum{\displaystyle\sum}

\newcommand{\pf}{\noindent{\it Proof \ }}
\newcommand{\epf}{$\quad$\hfill
\raisebox{0.11truecm}{\fbox{}}\par\vskip0.4truecm}

\begin{document}
\title[Affine Weyl groups and Frobenius manifolds]{Extended affine Weyl groups and Frobenius manifolds -- II}
\author{
{Boris Dubrovin \ \ Youjin Zhang \ \ Dafeng Zuo}}
\address{Dubrovin, SISSA, Via Beirut 2-4, 34014 Trieste, Italy}
\address{Zhang, Department of Mathematical Sciences, Tsinghua
University, Beijing 100084, P.R.China}
\address{Zuo, Department of Mathematical Sciences, Tsinghua
University, Beijing 100084, P.R.China and Department of
Mathematics, University of Science and Technology, Hefei 230026, P.R.China}

\email{dubrovin@sissa.it, yzhang@math.tsinghua.edu.cn,
\newline\mbox {\hskip2.97cm dfzuo@ustc.edu.cn}}
\date{}


\dedicatory{}

\begin{abstract}
For the root system of type $B_l$ and $C_l$, we generalize the result of \cite{DZ1998} by showing
the existence of a Frobenius manifold structure on the orbit space of the extended
affine Weyl group that corresponds to any vertex of the Dynkin diagram instead of a particular choice of \cite{DZ1998}.
\end{abstract}
\maketitle

\section{Introduction}
For an irreducible reduced root system $R$ defined in
$l$-dimensional Euclidean space $V$ with Euclidean inner product
$(~,~)$, we fix a basis of simple roots $\al_1,\dots,\al_l$ and
denote by  $\alpha_j^{\vee},\ j=1,2,\cdots,l$ the corresponding
coroots. The Weyl group $W$ is generated by the reflections
\beq
{\bx}\mapsto \bx-(\al_j^\vee,\bx)\al_j,\quad \forall\, {\bx}\in V,\ j=1,\dots,l.
\eeq The
semi-direct product of $W$ by the lattice of coroots yields the
affine Weyl group $W_a$ that acts on $V$ by the affine
transformations \beq \bx\mapsto w(\bx)+\sum_{j=1}^l m_j\av_j,\quad
w\in W,\ m_j\in \mathbb{Z}. \eeq
We denote by $\omega_1,\dots,\omega_l$ the fundamental weights that are defined by relations
\beq
(\omega_i,\al_j^\vee)=\delta_{ij},\quad i, j=1,\dots,l.
\eeq
Fixing a simple root $\al_k$,  we define an
{\it extended affine Weyl group} $\tw=\tw^{(k)}(R)$  as in
\cite{DZ1998}.
It acts on the extended space
$$
\tv=V\oplus \mathbb{R}
$$
and is generated by the transformations
\begin{equation}
x=(\bx,x_{l+1})\mapsto (w(\bx)+\sum_{j=1}^l m_j\av_j, \
x_{l+1}),\quad w\in W,\ m_j\in \mathbb Z,
\end{equation}
and
\begin{equation}
x=(\bx,x_{l+1})\mapsto (\bx+\w_k,\ x_{l+1}-1).
\end{equation}

Let us introduce coordinates $x_1,\dots,x_l$ on the space $V$ by
\beq
\bx=x_1\, \al_1^\vee+\dots+x_l\, \al_l^\vee.
\eeq
Denote by $f$ the determinant of the Cartan matrix of the root system $R$.
\begin{de}[\cite{DZ1998}] ${\mathcal A}={\mathcal A}^{(k)}(R)$ is
the ring of all $\ \tw$-invariant Fourier polynomials of the form
$$
\sum_{m_1,\dots,m_{l+1}\in\, \mathbb{Z}} a_{m_1,\dots,m_{l+1}} e^{2\pi i
(m_1 x_1+\cdots+m_l x_l+\frac1{f}\,m_{l+1} x_{l+1})}
$$
that are bounded in the limit
\begin{equation}
\bx=\bx^{0}-i\ \w_k\tau,\quad x_{l+1}=x_{l+1}^{0}+i\ \tau,\quad
\tau\to +\infty
\end{equation}
for any $\ x^{0}=(\bx^{0},x_{l+1}^{0})$.
\end{de}

For the fixed simple root $\al_k$, we introduce a set of numbers
\beq\label{zh-16}
d_j=(\omega_j, \omega_k),\quad j=1,\dots,l
\eeq
and define the following Fourier polynomials \cite{DZ1998}
\eqa
&&\ty_j(x)=e^{2\pi i d_j x_{l+1}}\, y_j(\bx),\quad j=1,\dots,l,\label{ip-a}\\
&&\ty_{l+1}(x)=e^{2\pi i x_{l+1}}.\label{ip-b}
\eeqa
Here
$y_1(\bx),\dots,y_l(\bx)$ are the $W_a$-invariant Fourier
polynomials defined by
\beq y_j(\bx)=\frac1{n_j} \sum_{w\in W}
e^{2\pi i (\omega_j, w(\bx))},\quad n_j=\#\{w\in W|e^{2\pi i
(\omega_j,w(\bx))}=e^{2\pi i (\omega_j,\bx)}\}.\label{fip}
\eeq
It was shown in \cite{DZ1998} that for some particular choices of the
simple root $\al_k$, a Chevalley type theorem holds true for the
ring ${\mathcal A}$, i.e., it is generated by
$\ty_1,\dots,\ty_{l+1}$, and thus the orbit space defined as
$\mathcal{M}={\rm Spec}{\mathcal A}$ of the extended affine Weyl
group $\widetilde{W}$ is an affine algebraic variety of dimension
$l+1$. {}Furthermore, in \cite{DZ1998} it was proved that on such an orbit space there
exists a {\it Frobenius manifold structure} whose potential is a
polynomial of $t^1,\dots, t^{l}, e^{t^{l+1}}$. Here $t^1,\dots,
t^{l+1}$ are the flat coordinates of the Frobenius manifold. For
the root system of type $A_l$, there is in fact no restrictions on
the choice of $\al_k$. However, for the root systems of
type $B_l, C_l, D_l, E_6, E_7, E_8, F_4, G_2$ there is
only one choice for each. Recall that the geometric structures on the orbit spaces
$\mathcal{M}$ generalize those that live on the orbit spaces of the finite Coxeter groups
discovered in \cite{SYS, S} and \cite{Du1}.

In \cite{slodowy} P. Solodowy showed that the Chevalley type theorem of \cite{DZ1998} can also be derived from
the results of
K. Wirthm\"uller \cite{wir}, and in fact it holds true for any choice of the base element $\al_k$, or equivalently,
for any fixed vertex of the Dynkin diagram. So we have

\begin{thm}[\cite{slodowy,wir}]\label{thm1}
The ring $\mathcal{A}$ is isomorphic to the ring of polynomials of
$\tilde{y}_1(x)$, $\cdots$, $\tilde{y}_{l+1}(x)$.  \end{thm}

A natural question, as it was pointed out in \cite{DZ1998, slodowy}, is
whether the geometric structures
revealed in \cite{DZ1998}
also exist on the orbit spaces of the extended affine Weyl groups for an arbitrary choice of $\al_k$?
The purpose of the present paper is to give an affirmative answer to this question for the root systems
of type $B_l, C_l$. We will show that
on the corresponding orbit spaces there also exist  Frobenius manifold structures with potentials that are polynomials
in $t^1,\dots, t^{l-1}, t^l, \frac1{t^l}, e^{t^{l+1}}$. Here $t^1,\dots, t^{l+1}$ are the flat coordinates
of the resulting Frobenius manifold.

The paper is organized as follows: in Sec.\ref{sec-2}, we
construct a flat pencil of metrics on each orbit space of the
extended affine Weyl group of the root system of type $C_l$ for any fixed vertex of the Dynkin diagram,
then in Sec.\ref{sec-3}
we study properties of the flat coordinates of the flat metric $(\eta^{ij})$.
In Sec.\ref{sec-4}, we prove the existence of a Frobenius manifold
structure on each orbit space. In Sec.\ref{sec-5} we give some examples.
In Sec.\ref{sec-6} we show that to the root system of type $B_l$ we can apply
a similar construction as the one for the root system of type $C_l$. The resulting Frobenius
manifolds are isomorphic to those that are obtained from $C_l$. Some concluding remarks are
given in the last section.

\noindent {\bf Acknowledgments.} The researches of B.D. were partially supported by European
Science Foundation Programme ``Methods of Integrable Systems, Geometry, Applied Mathematics"
(MISGAM), by the Marie Curie RTN ``European Network In Geometry, Mathematical physics and Applications'' (ENIGMA)
and by Italian Ministry of Universities and Scientific Researches research grant Prin2004
``Geometric methods in the theory of nonlinear waves and their applications''.
The researches of Y.Z. were partially supported by the Chinese National Science Fund
for Distinguished Young Scholars grant No.10025101 and the special Funds of Chinese Major Basic
Research Project ``Nonlinear Sciences". Y.Z. and D.Z. thank Abdus Salam International Center
for Theoretical Physics and SISSA where part of their work was done for the hospitality.
The authors thank Si-Qi Liu for his help on the proof of lemma \ref{lem3.1}.

\section{Flat pencils of metrics on the orbit spaces of ${\widetilde{W}}^{(k)}(C_l)$}\label{sec-2}

Let $\mathcal{M}$ be the orbit space defined as ${\rm Spec}{\mathcal A}$ of the extended affine Weyl
group $\widetilde{W}^{(k)}(C_l)$ for any fixed $1\le k\le l$. We choose the standard base $\al_1,\dots, \al_l$
of simple roots for the root system $C_l$ as given in \cite{bb}.
As in \cite{DZ1998} we define an indefinite  metric $(~,~)^{\sptilde}$
on $\widetilde V=V\oplus\mathbb R$ such that $\widetilde V$ is the orthogonal direct sum of $V$ and $\mathbb R$,
$V$ endowed with the $W$-invariant Euclidean metric
\beq
(d x_m, dx_n)^{\sptilde}=\frac{m}{4\pi^2},\quad 1\le m\le n\le l
\eeq
and $\mathbb R$ endowed with the metric
\beq
(dx_{l+1},dx_{l+1})^{\sptilde}=-\frac1{4\pi^2 d_k}=-\frac1{4 k \,\pi^2}.
\eeq
Here the numbers $d_j$ are defined in (\ref{zh-16}) and take the values
\beq
d_1=1,\dots, d_{k-1}=k-1,\ d_j=k,\quad k\le j\le l. \eeq
The $W_a$-invariant Fourier
polynomials $y_1(\bx),\dots, y_l(\bx)$ that are defined in (\ref{fip}) have the expressions
\beq\label{zh11}
y_j({\bx})=\sigma_j(\xi_1,\cdots,\xi_l),
\eeq
where
\beq
\xi_{{j}}={e^{2\,i\pi\, \left( x_{{j}}-x_{{j-1}} \right)
}}+{e^{-2\,i \pi\, \left( x_{{j}}-x_{{j-1}} \right) }},~~
x_{0}=0, ~~j=1,\cdots,l
\eeq
and $\sigma_j(\xi_1,\dots,\xi_l)$ is the $j$-th elementary
symmetric polynomial of $\xi_1,\cdots,\xi_l$.
For the reason that will be clear later, we will use in what follows the following
set of generators for the the ring of $W_a$-invariant Fourier
polynomials which we still denote by $y_1,\dots, y_l$:
\beq\label{zh10}
{y}_j({\bx})=\sigma_j(\zeta_1,\cdots,\zeta_l)
\eeq
with
\beq\label{zh7b}
\zeta_j=\xi_j+2,\quad j=1,\dots,l.
\eeq
Consequently we have a set of generators for the ring ${\mathcal A}={\mathcal A}^{(k)}(C_l)$
\eqa
&&\ty_j(x)=e^{2\pi i d_j x_{l+1}}\, y_j(\bx),\quad j=1,\dots,l,\label{zh8-a}\\
&&\ty_{l+1}(x)=e^{2\pi i x_{l+1}},\label{zh8-b} \eeqa defined in
the same way as in (\ref{ip-a}), (\ref{ip-b}).
They form
a global coordinates on ${\mathcal M}$. As in \cite{DZ1998}, we
introduce the following local coordinates on ${\mathcal M}$:
\beq\label{zh9} y^1=\ty_1,\dots, y^l=\ty_l,\ y^{l+1}=\log
\ty_{l+1}=2 \pi i\, x_{l+1}. \eeq They live on a covering
$\widetilde{\mathcal M}$ of ${\mathcal
M}\setminus\{\ty_{l+1}=0\}$. The projection \beq P: \widetilde
V\to \widetilde{\mathcal M} \eeq induces a symmetric bilinear form
on $T^{*}\widetilde{\mathcal{M}}$
\begin{equation}
(d y^i,d y^j)^{\sptilde}\equiv
g^{ij}(y):=\dsum_{a,b=1}^{l+1}\dfrac{\p y^i}{\p x^a}\dfrac{\p
y^j}{\p x^b}(dx^a,dx^b)^{\sptilde}.
\end{equation}
Denote
\beq
\Sigma=\{y|\det(g^{ij}(y))=0\}.
\eeq
It turns out that $\Sigma$ is an analogue of the discriminant. Namely, as
it was shown in \cite{DZ1998}, $\Sigma$ is the $P$-image of the union of hyperplanes
\begin{equation}
\{(\bx,x_{l+1})|(\beta,\bx)=m\in\mathbb{Z},\ x_{l+1}=\text {arbitrary}\},
\quad \beta\in \Phi^+,\tag{2.3}
\end{equation}
where $\Phi^+$ is the set of all positive roots.

\begin{lem} The functions $g^{ij}(y)$ are weighted homogeneous polynomials in
${y}^1,\cdots$,~$y^l$, $e^{y^{l+1}}$ of the degree
\begin{equation}\deg g^{ij}=\deg y^i+\deg y^j,\end{equation}
where $\deg y^j=d_j$ and $\deg y^{l+1}=d_{l+1}=0$.
\end{lem}
\pf The assertion follows immediately from the Theorem \ref{thm1} which says that the
functions $\ty_1,\dots,\ty_{l+1}$ form a set of generators of the
ring ${\mathcal A}$. \epf

From this lemma we see that $\Sigma$ is an algebraic subvariety in $\mathcal M$ and the
matrix $(g^{ij})$ is invertible in ${\mathcal M}\setminus \Sigma$. The inverse matrix
$(g^{ij})^{-1}$ defines a flat metric on ${\mathcal M}\setminus \Sigma$.

Let us introduce the following new coordinates on ${\mathcal M}$
\beq
\ta^j=\left\{\begin{array}{ll}e^{k\, y^{l+1}},&j=0,\\
y^j e^{(k-j) y^{l+1}}, &j=1,\cdots,k-1,\\
y^j,  &j=k,\cdots,l
\end{array}\right. .
\eeq
Denote
\beq \om_j=2\pi
i(x_j-x_{j-1}),\quad \om_{l+1}=y^{l+1}, \quad j=1,\cdots,l.
\eeq
In the coordinates $\om_1,\dots, \om_{l+1}$ the indefinite metric on
$\widetilde V$ has the form
\beq \left((d \om_i,d \om_j)^{\sptilde} \right)=\hbox{diag}(-1,\dots,-1,\frac1{k}). \eeq Define
\begin{equation}
P(u):=\dsum_{j=0}^{l} u^{l-j}\theta^j=e^{k \om_{l+1}}
\prod_{j=1}^{l}(u+\zeta_j).\label{P1}
\end{equation}
We can easily verify that the function $P(u)$ satisfies \eqa
&&\dfrac{\p P(u)}{\p
\om_a}=\dfrac{1}{u+\zeta_a}P(u)(e^{\om_a}-e^{-\om_a}),\quad
1\leq a \leq l;\label{eqp-a}\\
&&\dfrac{\p P(u)}{\p \om_{l+1}}=k P(u), \quad P'(u):=\dfrac{\p
P(u)}{\p u}=P(u)\dsum_{a=1}^{l}\frac{1}{u+\zeta_a}.\label{eqp-b}
\eeqa By using these identities, we have
\begin{lem} The following formulae hold true for the generating functions of the
metric $(g^{ij})$ and the contravariant components of its
Levi-Civita connection in the coordinates $\ta^0,\dots,
\ta^{l}$
\eqa
&&\dsum_{i,j=0}^{l}(d\ta^i,d\ta^j)^{\sptilde}u^{l-i}v^{l-j}=(d P(u), d P(v))^{\sptilde}\nn\\
&&\quad
=(k-l)P(u)P(v)+\dfrac{u^2+4u}{u-v}P'(u)P(v)-\dfrac{v^2+4v}{u-v}P(u)P'(v),\nn
\eeqa
\vskip -0.75truecm
\eqa &&\dsum_{i,j,m=0}^{l}\Gamma^{ij}_m(\ta)
d\ta^m u^{l-i}v^{l-j}=\dsum_{a,b,m=1}^{l+1}\dfrac{\p P(u)}{\p
\om_a}\dfrac{\p^2 P(v)}{\p \om_b \p
\om_m}d\om_m(d\om_a,d\om_b)\nn\\
&&\quad =(k-l)P(u)dP(v)+\dfrac{u^2+4u}{u-v}P'(u)dP(v)-\dfrac{v^2+4v}{u-v}P(u)dP'(v)\nn\\
&&\qquad+\dfrac{2u+uv+2v}{(u-v)^2}P(v)dP(u)-\dfrac{2u+uv+2v}{(u-v)^2}P(u)dP(v).\nn
\eeqa
Here
$\Gamma^{ij}_m(\theta)=-\dsum_{s=1}^{l+1}g^{is}(\theta) \Gamma^j_{sm}(\theta)$.
\end{lem}
\pf  By using \eqref{eqp-a}, \eqref{eqp-b} we have
\begin{equation*}
\begin{array}{rl}
&(dP(u),dP(v))^{\sptilde}=\dfrac{1}{k}\dfrac{\p P(u)}{\p \om_{l+1}}\dfrac{\p P(v)}{\p
\om_{l+1}}-\dsum_{a=1}^l\dfrac{\p P(u)}{\p \om_a}\dfrac{\p
P(v)}{\p \om_a}\\
&=k P(u)P(v)-\dsum_{a=1}^l P(u)P(v)\frac{\zeta_a^2-4\zeta_a}{(u+\zeta_a)(v+\zeta_a)}\\
&=\dsum_{s=1}^l P(u)P(v)\frac{v}{u-v}(1-\frac{v}{v+\zeta_a})-\dsum_{a=1}^l P(u)P(v)
\frac{u}{u-v}(1-\frac{u}{u+\zeta_a})\\
&\quad+\dsum_{a=1}^l P(u)P(v)\frac{1}{u+\zeta_a}\frac{4u}{u-v}-\dsum_{a=1}^l
P(u)P(v)\frac{1}{v+\zeta_a}\frac{4v}{u-v}+kP(u)P(v)\\
&=(k-l)P(u)P(v)+\dfrac{u^2+4u}{u-v}P'(u)P(v)-\dfrac{v^2+4v}{u-v}P(u)P'(v).
\end{array}
\end{equation*}
So we proved the first formula, the second formula can be proved
in the same way. The lemma is proved. \epf The above lemma shows
that in the coordinates $\ta^0,\dots, \ta^{l}$ the functions
$g^{ij}(\ta)$ are quadratic polynomials, and the contravariant
components $\Gamma^{ij}_m$ are homogeneous linear
functions\footnote{These metrics give rise to a quadratic Poisson
structure on the space of ``loops" $\{S^1\to M\}$ (see \cite{Du2} for the details):
\beq
\{\theta^i(s), \theta^j(s')\}=g^{ij}(\theta(s))\delta'(s-s')
+\Gamma_m^{ij}(\theta(s))
\theta_s^m\delta(s-s').\nn
\eeq
We plan to study such important  class of quadratic metrics and Poisson structures
in a separate publication.}. It
reveals the following important properties of the flat metric:

\begin{cor}
In the coordinates $y^1,\dots, y^{l+1}$ the functions
$\Gamma_m^{ij}(y)$ are weighted homogeneous polynomials of degree
\begin{equation}
\deg \Gamma_m^{ij}(y)=d_i+d_j-d_m.
\end{equation}
\end{cor}

\begin{cor}\label{lem2.1}
In the coordinates $y^1,\dots, y^{l+1}$ the polynomials $g^{ij}(y)$ and $\Gamma^{ij}_m(y)$
are at most linear in $y^k$.
\end{cor}

Now let us define a symmetric bilinear form on $T^{*}\mathcal{M}$ by
\begin{equation}
<d y^i, dy^j>:=\eta^{ij}(y)=\dfrac{\p g^{ij}(y)}{\p y^k}.\end{equation}
\begin{lem}\label{lem2.2}
The matrix $(\eta^{ij}(y))$ has the form
\beq
\left(\begin{array}{cccccccccccccc}
0&0&0&\cdots&0& k & P_{1}&0&0&\cdots&0&0\\
0&0&0&\cdots&k&R_{1}& P_{2}&0&0&\cdots&0&0\\
0&0&0&\cdots&R_1&R_{2}&P_{3}&0&0&\cdots&0&0\\
\vdots&\vdots&\vdots&\vdots&\vdots&\ddots&\vdots&\vdots&\vdots&\ddots&\vdots&\vdots\\
k&R_{1}&R_2&\cdots&{}&R_{k-2}&P_{k-1}&0&0&\cdots&0&0\\
P_{1}&P_{2}&P_3&\cdots&P_{k-2}&P_{k-1}&P_{k}&0&0&\cdots&0&1\\
0&0&0&\cdots&0&0&0&Q_1&Q_2&\cdots&Q_{l-k}&0\\
0&0&0&\cdots&0&0&0&Q_2&Q_3&\cdots&0&0\\
\vdots&\vdots&\vdots&\ddots&\vdots&\vdots&\vdots&\vdots&\vdots&\cdots&\vdots&\vdots\\
0&0&0&\cdots&0&0&0&Q_{l-k}&0&\cdots&0&0\\
0&0&0&\cdots&0&0&1&0&0&\cdots&0&0
\end{array}\right).\label{zh1}
\eeq
Here
\eqa
&&R_{j}=4(k-j+1) y^{j-1} e^{y^{l+1}} +(k-j) y^j,\nn\\
&&P_j=4(k-j+1) y^{j-1} e^{y^{l+1}},\nn\\
&&Q_m=4 m y^{k+m}+(1-\delta_{m,l-k}) (m+1) y^{k+m+1},\nn\\
&&1\le j\le k,\quad 1\le m\le l-k\nn
\eeqa
and we assume $y^0=1$.
\end{lem}

It follows from the above lemma that
  $$\det (\eta^{ij})=(-1)^l k^{k-1} 4^{l-k} (l-k)^{l-k} (y^l)^{l-k},$$
so $(\eta^{ij})^{-1}$ defines a metric on ${\mathcal M}\setminus \{y\in{\mathcal M}|y^l=0\}$.

\begin{thm}The space $\mathcal{M}$ carries a flat
pencil of metrics (bilinear forms on $T^{*}M$)
$$
g^{ij}(y) \ \ {  {and}} \ \ \eta^{ij}(y)=\dfrac{\p g^{ij}(y)}{\p y^k},
$$
i.e., any linear combination $g^{ij}+\lambda \,\eta^{ij}$ defines a flat metric on certain open subset
of ${\mathcal M}$ and the contravariant components of the Levi-Civita connection of $(\eta^{ij})$ is given by
\beq
\gamma_m^{ij}(y)=\dfrac{\p \Gamma_m^{ij}(y)}{\p y^k}.
\eeq
The metric $(g^{ij}(y))$ does not degenerate on ${\mathcal M}\setminus \Sigma$ and
the metric $(\eta^{ij}(y))$ does not degenerate on ${\mathcal M}\setminus \{y\in{\mathcal M}|y^l=0\}$.
\end{thm}
\pf The result follows, applying
Lemma D.1 of \cite{Du1}, from the fact that in the coordinates $y^1,\dots, y^{l+1}$ the flat metric $g^{ij}$
and the contravariant components of its Levi-Civita connection
depend at most linearly on $y^k$ . The theorem is proved. \epf

\begin{rem}
Our particular choice of the basis of the ${\widetilde W}^{(k)}(C_l)$-invariant Fourier polynomials
(\ref{zh10})--(\ref{zh8-b})
ensures that the components of the flat metric $(g^{ij}(y))$  are at most linear in $y^k$. This linearity is the most
crucial step in the construction of the above flat pencil of metrics.
If we choose the basis of the ${\widetilde W}^{(k)}(C_l)$-invariant Fourier polynomials by using \eqref{zh11},
\eqref{zh8-a} and \eqref{zh8-b}, then we lose this linearity property of the functions $(g^{ij}(y))$ and the
construction of the flat metric $(\eta^{ij}(y))$ becomes obscure.
\end{rem}


\section{Flat coordinates of the metric $(\eta^{ij})$}
\label{sec-3}
In this section, we will show that the flat coordinates of the metric
$(\eta^{ij})$ are algebraic functions of
$y^1,\dots, y^{l+1}, e^{y^{l+1}}$. To this end, we first perform changes of coordinates to simplify
the matrix $(\eta^{ij})$.

\begin{lem}\label{lem3.1}
There exists a system of new coordinates $z^1,\dots, z^{l+1}$ of
the form \eqa &&z^j=y^j+p_j(y^1,\dots,y^{j-1},e^{y^l+1}), ~~1\leq j \leq k, ~~
z^{l+1}=y^{l+1},
\label{zh2-a}\\
&& z^j=y^j+\sum_{m=j+1}^l c^j_m y^m,\quad k+1\le j\le
l,\label{zh2-b} \eeqa where $p_j$ are homogeneous polynomials of
degree $d_j$ and $c^j_m$ are some constants such that in the
new coordinates  the metric $(\eta^{ij})$ still has the form
(\ref{zh1}) with the entries replaced by \beq R_{j}=0,\quad
P_j=0,\quad Q_m=4 m z^{k+m},\quad 1\le j\le k,\quad 1\le m\le
l-k.\label{zh3} \eeq
\end{lem}
\pf Let us first note that the $(k+1)\times (k+1)$ matrix $({\tilde{\eta}}^{ij})$
which has elements
\beq
{\tilde {\eta}}^{ij}=\eta^{ij}(y),\ {\tilde{\eta}}^{k+1,m}={\tilde{\eta}}^{m,k+1}=\delta_{j,k},\quad
1\le i,j\le k,\ 1\le m\le k+1 \nn
\eeq
coincides, under renaming of the coordinate $y^{l+1}\mapsto y^{k+1}$,
 with the  matrix $(\eta^{ij}(y))_{(k+1)\times (k+1)}$ that is constructed as in the last section
with respect to the extended affine Weyl group ${\widetilde
W}^{(k)}(C_k)$. Thus by using the results of \cite{DZ1998} we can
find homogeneous polynomials $p_j, 1\le j\le k$ such that under
the change of coordinates \eqref{zh2-a} and $z^j=y^j$, $k+1\le j\le
l$ the matrix $(\eta^{ij}(z)$ has the form \eqref{zh1} with
entries \eqa
&&R_{j}=0,\quad P_j=0,\quad Q_m=4 m z^{k+m}+(1-\delta_{m,l-k}) (m+1) z^{k+m+1},\nn\\
&&\quad 1\le j\le k,\quad 1\le m\le l-k.\nn
\eeqa

To finish the proof of the lemma, we need to perform a second change of coordinates.
To this end, denote by $\Psi$ a $n\times n$ matrix with entries
 as linear functions of $a^1,\dots, a^n$
\begin{equation}
\psi^{ij}(a)=4(i+j-1)a^{i+j-1}+(i+j)a^{i+j},\quad i,j\ge 1.\nn
\end{equation}
Here $a^m=0 \ {\rm for}\  m\ge n+1$.
We are to find a linear transformation of the triangular form
\beq
a^j=\dsum_{m=j}^n B_{m}^jb^m,\quad B^j_j=1,~~j\ge 1 \nn
\eeq
such that
\eqa
&&\dsum_{r,s=1}^n 4(r+s-1)b^{r+s-1}\dfrac{\p a^i}{\p b^r}\dfrac{\p
a^j}{\p b^s}\nn\\
&&=4(i+j-1) \sum_{m=i+j-1}^n B^{i+j-1}_m b^m
+(i+j) \sum_{m={i+j}}^n B^{i+j}_m b^m.\nn
\eeqa
Equivalently, the constants $B^i_j$ must satisfy the relations
\eqa
&&4(i+j-1)B^{i+j-1}_m+(i+j)B^{i+j}_m=4m
\dsum_{\alpha+\beta=m+1}B^i_\alpha B^j_\beta,\nn\\
&&\quad i+j\le m\le n.\label{TR8}
\eeqa
Introduce the generating functions
\beq\label{TR9}
f^i(t)=\sum_{\alpha\ge0}B^i_{i+\alpha}t^\alpha,\quad i=1,2,\dots.
\eeq
Then the relations in \eqref{TR8} can be encoded into the following equations of $f^i(t)$:
\beq\label{TR10}
4(i+j-1)t^{i+j-2}f^{i+j-1}+({i+j})
t^{i+j-1}f^{i+j}= 4\frac{d}{dt}\left(t^{i+j-1}f^if^j\right).
\eeq
This system of equations has the following solution which was obtain by Si-Qi Liu
\beq\label{TR11}
f^i(t)=\cosh\left(\frac{\sqrt{t}}2\right)\left(\frac{2\sinh\left(\frac{\sqrt{t}}2\right)}{\sqrt{t}}\right)^{2i-1}.
\eeq
From the above result we derive the existence of constants $c^j_m,\ k+1\le j\le l,\ j+1\le m\le l$
such that under the change of coordinates
\beq
z^i\mapsto z^i,\ i=1,\dots, k, l+1,\quad
z^j\mapsto z^j+\sum_{m=j+1}^l c^j_m z^m,\quad k+1\le j\le l,\nn
\eeq
the matrix $(\eta^{ij}(z)$ has the form (\ref{zh1}) and with entries given by (\ref{zh3}).
The lemma is proved. \epf

\begin{lem}
Under the change of coordinates
\eqa
&&w^i=z^i, \quad i=1,\dots,k,\ l+1,\\
&&w^{k+1}=z^{k+1} (z^l)^{-\frac1{2 (l-k)}},\
w^j=z^j(z^{l})^{-\frac{j-k}{l-k}}, \
w^{l}=(z^{l})^{\frac{1}{2(l-k)}},\\
&&j=k+2\cdots,l-1,\nn
\eeqa
the components of the metric $(\eta^{ij}(z))$ are transformed to the form
$$
{\rm {diag}}(A_{(k-1)\times (k-1)}, B_{(l-k+2)\times (l-k+2)}),
$$
where the matrix $A$ has entries $A^{ij}=\delta_{i,k-i} k$ and the upper triangular
matrix $B$ has the form
\beq
B=\left(\begin{array}{cccccccccccc}
0&0&0& 0&0&\cdots &0 &0&1 \\
0&0&0& 0&0&\cdots &0 &2& \\
0&0&S_{k+3}&S_{k+4}&\cdots&S_{l-1}&S_{l}&&\\
0&0&S_{k+4}&S_{k+5}&\cdots &S_l&&\\
\vdots&\vdots &\vdots&&&&&&\\
0&0&S_l&& &&&&\\
0&2&&&&&&&\\
1&&&&&&&&
\end{array}\right).
\label{TR15}\eeq
with
\beq
S_{k+j}=4 j (w^l)^{-2} w^{k+j},~~ S_l=4 (l-k) {(w^l)}^{-2},\quad 3\le j\le l-k-1.
\eeq
\end{lem}
\pf By a straightforward calculation. \epf

\begin{prop}
In the coordinates $w^1,\dots,w^{l+1}$ the Christoffel symbols of the metric $(\eta^{ij})$ have the
following properties
\begin{enumerate}
\item \ $\gamma^{m}_{ij}=0$ \ for\ $m=1,\dots, k, l, l+1,\
i,j=1,\dots,l+1;$
\item \ $\gamma^{k+1}_{ij}=-\frac{\pal\eta_{ij}(w)}{\pal w^l}$ are
weighted homogenous polynomials in $w^{k+3},\dots, w^l;$
\item \ $\gamma^m_{ij}$ \ are weighted homogenous polynomials in\
$w^{k+3},\dots, w^l$ \ for\ $k+2\le m\le l-1$ and
$i,j=1,\dots, l-1, l+1;$
\item \ $\gamma^m_{lj}=\frac1{w^l}
\delta^m_j$, for \ $k+2\le m\le l-1,\ 1\le j\le l+1$.
\end{enumerate}
\end{prop}
\pf The first three properties of $\gamma^m_{ij}$ follow easily from the
simple form of the matrix $(\eta^{ij}(w))$. To prove the last property, we only need to
note that
\beq
\gamma^m_{lj}=\frac12 \sum_{s=1}^{l+1} \eta^{ms}\frac{\pal\eta_{sj}}{\pal w^l}
=\frac12 \sum_{s=1}^{l+1} \frac{2}{w^l} \eta^{ms} \eta_{sj}=\frac1{w^l} \delta^m_j.
\eeq
The proposition is proved. \epf

\begin{thm}
We can choose the flat coordinates of the metric $(\eta^{ij}(w)$
in the form
\eqa &&t^1=w^1,\dots, t^k=w^k,\ t^l=w^l,\
t^{l+1}=w^{l+1},\\&& t^{k+1}=w^{k+1}+w^l\,h_{k+1}(w^{k+2},\dots,w^{l-1}),\\
&&
t^j=w^l (w^j+h_j(w^{j+1},\dots,w^{l-1})).
\eeqa Here $h_j$ are weighted homogeneous
polynomials of degree $\frac{k\,(l-j)}{l-k}$  for $j=k+1,\dots,l-2$ and $h_{l-1}=0$.
\end{thm}
\pf To find the flat coordinates $t=t(w)$, we need to solve the following system of PDEs
\beq\label{zh5}
\frac{\pal^2 t}{\pal w^i\pal w^j}-\sum_{m=1}^{l+1}\gamma^m_{ij}
   \frac{\pal t}{\pal w^m}=0,\quad i,j=1,\dots,l+1.
\eeq
From the above proposition we easily see that $t^1,\dots, t^k, t^{l+1}$ are $k+1$ solutions of the above
system. We still need to find $l-k$ independent solutions $t^{k+1},\dots, t^l$. Introduce the $(l-k)\times (l-k)$
matrix
\beq
\Phi=(\phi^i_j),\quad \phi^i_j=\frac{\pal t^{k+i}}{\pal w^{k+j}}.
\eeq
Then the system (\ref{zh5}) is reduced to
\beq\label{zh6}
\pal_m\Phi=\Phi A_m,\quad \pal_m=\frac{\pal}{\pal w^m}, \quad m=k+1,\dots,l.
\eeq
They are regular at ${\bf w}=(w^{k+1},\dots,w^l)=0$ except the system with
$m=l$. In this case the coefficient
matrix has the simple form
\beq
A_l={\mbox{diag}}(0,\frac1{w^l},\dots,\frac1{w^l},0).
\eeq
Now assume that $\Phi$ has the form
\beq
\Phi=\Psi \,{\mbox{diag}}(1,w^l,\dots,w^l,1).
\eeq
Then the systems in (\ref{zh6}) are converted to
\beq
\pal_m\Psi=\Psi B_m,\quad \pal_l\Psi=0,\quad m=k+1,\dots,l-1.
\eeq
The entries of the coefficient matrices $B_m$ are now weighted homogeneous
polynomials
of $w^{k+1},\dots, w^l$. Thus we can find a unique analytic at ${\bf w}=0$ solution $\Psi$ of the above systems
such that
\beq
\left.\Psi\right|_{{\bf w}=0}={\mbox{diag}}(1,\dots,1).
\eeq
From the weighted homogeneity of the coefficient matrices $B_m$ it follows that
the elements of $\Psi$ are also weighted homogeneous. Since $\deg w^j>0$ for $j=k+1,\dots, l$
we see that they are in fact polynomials of $w^{k+1},\dots, w^l$. Thus the result of the theorem follows.
The theorem is proved. \epf

Due to the above construction, we can associate the following
natural degrees to the flat coordinates
\eqa
&&\tilde d_j=\deg t^j:=\frac{j}{k},\quad 1\le j\le k,\label{zh18a}\\
&& \tilde d_m=\deg t^m:=\dfrac{2l-2m+1}{2(l-k)},\quad  k+1\le m\le l,\label{zh18b}\\
&&\tilde d_{l+1}=\deg t^{l+1}:=0,\label{zh18c}
\eeqa
and we readily have the
following corollary
\begin{cor}\label{cor3.5}
In the flat coordinates $t^1,\dots, t^{l+1}$, the entries of
the matrix $(\eta^{ij})$ has the form \beq \label{zh7}
\eta^{ij}=\left\{\begin{array}{lll} k,\quad & j=k-i,\ & 1\le i\le k-1,\\
4 (l-k),\quad  &j=k+l+1-i,\ &k+2\le i\le l-1,\\
1,&i=l+1, j=k\ & {\mbox or}\ i=k,\ j=l+1,\\
2,&i=l, j=k+1\ &{\mbox or}\ i=k+1,\ j=l.
\end{array}\right.\eeq
 The entries of the matrix $(g^{ij}(t))$ and the
Christoffel symbols $\Gamma^{ij}_m(t)$ are weighted homogeneous
polynomials of $t^1,\dots, t^l, \dfrac{1}{t^l}, e^{t^{l+1}}$ of
degrees $\td_i+\td_j$ and $\td_i+\td_j-\td_m$ respectively.
In particular,
\beq \label{zh8} \begin{array}{ll}g^{m,\,l+1}=\td_m
t^m, \quad & 1\leq m
\leq l,\quad g^{l+1, \,l+1}=\dfrac{1}{k},\\
\Gamma_{j}^{l+1,i}=\td_j\,\delta_{i,j},\quad & 1\leq i,j \leq l+1.
\end{array}
\eeq
\end{cor}

\begin{rem} For the orbit spaces of finite reflection groups flat coordinates were constructed by Saito,
Yano and Sekiguchi in \cite{SYS} (see also \cite{S}).
\end{rem}

The numbers $\td_1,\td_{l+1}$ satisfy a duality relation that is
similar to that of \cite{DZ1998}. To describe this duality
relation, let us delete the $k$-th vertex of the Dynkin diagram
${\mathcal R}$ and obtain two components ${\mathcal R}\setminus
{\al_k}={\mathcal R}_1\cup {\mathcal R}_2$. On each component we
have an involution \beq i\mapsto i^*, \quad i=1,\dots, k-1 \
{\mbox{and}}\ i=k+1,\dots,l \eeq given by the reflection with
respect to the center of the component. We also define \beq
k^*=l+1,\quad (l+1)^*=k, \eeq then we have \beq
\td_i+\td_{i^*}=1,\quad i=1,\dots, l+1, \eeq and from the above
corollary we see that $\eta^{ij}$ is a nonzero constant iff
$j=i^*$.
\section{The Frobenius manifold structure on the orbit space of $\widetilde{W}^{(k)}(C_l)$}
\label{sec-4}

Now we are ready to describe the Frobenius manifold structure on the orbit space of
the extended affine Weyl group ${\widetilde W}^{(k)}(C_l)$. Let us first recall the definition of Frobenius
manifold, see \cite{Du1} for details.

\begin{de} A {\it Frobenius algebra} is a pair $(A, <~,~>)$ where $A$
is a commutative associative algebra with a unity $e$ over a field $\mathcal{K}$ (in our case
$\mathcal{K}={\mathbb C}$) and $<~,~>$ is a $\mathcal{K}$-bilinear
symmetric nondegenerate {\it invariant} form on $A$, i.e.,
$$
<x\cdot y, z> = <x, y\cdot z>,\quad \forall\ x, y, z \in A.
$$

\end{de}
\begin{de}A Frobenius  structure of charge $d$
on an n-dimensional manifold $M$ is a structure of Frobenius algebra on the tangent spaces $T_tM=(A_t,<~,~>_t)$
depending (smoothly, analytically etc.) on the point $t$. This structure satisfies the following axioms:
 \begin{itemize}
\item[FM1.] The metric $<~,~>_t$ on $M$ is flat, and the unity vector field $e$ is covariantly constant, i.e.,
$\nabla e=0$. Here we
denote $\nabla$ the Levi-Civita connection for this flat metric.
\item[FM2.] Let $c$ be the 3-tensor $c(x,y,z):=<x\cdot y, z>$, $x,\, y,\,
z\in T_tM$. Then the 4-tensor $(\nabla_w c)(x,y,z)$ is symmetric in
$x,\, y,\, z, \, w \in T_tM$.
\item[FM3.] The existence on $M$ of a vector field $E$, called the Euler vector field,
which satisfies the conditions $\nabla\nabla E=0$ and
$$
[E, x\cdot y] -[E,x]\cdot y -x\cdot [E,y] = x\cdot y,
$$
$$
E<x,y>-<[E,x],y>-<x,[E,y]>=(2-d)<x,y>
$$
for any vector fields $x, y$ on $M$.
\end{itemize}
A manifold $M$ equipped with a Frobenius structure on it is called a Frobenius manifold.
 \end{de}

Let us choose locally flat coordinates $t^1,\cdots t^n$ for the
invariant flat metric, then locally there exists a function
$F(t^1,\cdots,t^n)$, called the {\em potential} of the Frobenius
manifold, such that \beq < u\cdot v,w>=u^i v^j w^s \frac{\p^3
F}{{\p t^i}{\p t^j}{\p t^s}} \label{WDVV0} \eeq for any three
vector fields $u=u^i\frac{\p}{\p t^i}$, $v=v^j\frac{\p}{\p t^j}$,
$w=w^s\frac{\p}{\p t^{s}}$. Here and in what follows summations
over repeated indices are assumed. By definition, we can also
choose the coordinates $t^1$ such that $ e=\frac{\pal}{\pal
t^1} $. Then in the flat coordinates the components of the
flat metric can be
expressed in the form
\beq \frac{\p^3 F}{{\p t^1}{\p t^i}{\p
t^j}}=\eta_{ij}=<\frac{\p}{\p t^i},\frac{\p}{\p t^j}>,\quad i,j=1,\dots, n. \label{WDVV1} \eeq
The associativity of the Frobenius algebras is equivalent to the
following overdetermined system of equations for the function $F$
\beq\label{WDVV2} \frac{\p^3 F}{{\p t^i}{\p t^j}{\p
t^\lambda}}\eta^{\lambda\mu} \frac{\p^3 F}{{\p t^\mu}{\p t^k}{\p
t^m}}=\frac{\p^3 F}{{\p t^k}{\p t^j}{\p
t^\lambda}}\eta^{\lambda\mu} \frac{\p^3 F}{{\p t^\mu}{\p t^i}{\p
t^m}} \eeq for arbitrary indices $i,j,k,m$ from $1$ to $n$.

In the flat coordinates the Euler vector field $E$ has the form
\beq E=\sum_{i=1}^n ({d}_j^i t^j+r_i)\frac{\pal}{\pal t^i} \eeq
for some constants ${d}_j^i, r_i,\, i=1,\dots,n$ which satisfy
$$
{d}_1^i=\delta_1^i,\quad r_1=0.
$$
From the axiom FM3 it follows that the potential $F$ satisfies
the quasi-homogeneity condition \beq \label{WDVV3}
\mathcal{L}_EF=(3-d)F+\frac12 A_{ij}\,t^i t^j+B_i\,t^i+\mbox{constant}.
\eeq
The system \eqref{WDVV1}--\eqref{WDVV3} is called the {\sl $WDVV$
equations of associativity} which is equivalent to the above
definition of Frobenius manifold in the chosen system of local
coordinates.

In our examples the constant matrix $d_i^j$ is always diagonal,
$d_i^j=\hat{d}_i\delta_i^j$.

Let us also recall an important geometrical structure on a Frobenius manifold $M$,
the {\sl intersection form} of $M$.
This is a symmetric
bilinear form $(~,~)^*$ on $T^*M$ defined by the formula
\beq (w_1,w_2)^*=i_E(w_1\cdot w_2),\eeq
here the product of two 1-forms $w_1$, $w_2$ at a point $t\in M$ is defined
by using the algebra structure on $T_tM$
and the isomorphism
\beq T_tM\to T_t^*M\eeq
 established by the invariant flat metric $<~,~>$. In the flat coordinates
  $t^1,\cdots,t^n$ of the invariant
 metric,  the intersection form can be represented by
 \beq
(dt^i,dt^j)^*=\mathcal{L}_EF^{ij}=(d+1-\hat{d}_i-\hat{d}_j)F^{ij}+A^{ij},
 \eeq
where
\beq
A^{ij}=\eta^{ii'}\eta^{jj'} A_{i'j'},\quad
F^{i j}=\eta^{i i'}\eta^{j j'}\dfrac{\p^2 F}{{\p
t^{i'}}{\p t^{j'}}}\eeq
and $F(t)$ is the potential of the Frobenius manifold.
Denote by $\Sigma\subset M$ the {\em discriminant} of $M$ on which the intersection form degenerates, then
an important property of the intersection form is that on $M\setminus\Sigma$ its inverse defines a new flat metric.

\begin{thm}There exists a unique Frobenius structure of charge $d=1$ on the orbit
space $\mathcal{M} \setminus\{t^l=0\}$ of ${\widetilde{W}}^{(k)}(C_l)$ polynomial in
$t^1,t^2, \cdots, t^l,  \dfrac{1}{t^{l}}, e^{t^{l+1}}$ such that
\begin{enumerate}
\item The unity vector field $e$ coincides with $\dfrac{\p}{\p
y^k}=\dfrac{\p}{\p t^k}$;
 \item The Euler vector field has the form
\begin{equation}E=\dsum_{\alpha=1}^{l}\td_\alpha t^\alpha \dfrac{\p}{\p t^\alpha}
+\dfrac{\p}{\p t^{l+1}}\label{zz1}
\end{equation}
where $\td_1,\dots,\td_{l}$ are defined in (\ref{zh18a})--(\ref{zh18c}).
\item The invariant flat metric and the intersection form of the Frobenius structure coincide respectively with the
metric $(\eta^{ij})$ and $(g^{ij}(t))$ on $\mathcal{M}\setminus\{t^l=0\}$.
\end{enumerate}
\end{thm}

\pf By following the lines of the proof of Lemma 2.6 that is given in \cite{DZ1998} we
can show the existence of a unique weighted homogeneous
polynomial
$$G:=G(t^1,\dots,t^{k-1},t^{k+1},\dots,t^l,\frac{1}{t^l},e^{t^{l+1}})$$
of degree $2$ such that the function \beq\label{zz2}
F=\frac{1}{2}(t^k)^2t^{l+1}+\frac{1}{2}t^k\dsum_{i,j\ne k}\eta_{i
j}\,t^i t^j+G \eeq satisfies the equations \beq \label{zz3} g^{i
j}=\mathcal{L}_E F^{i j},\quad \Gamma^{i j}_m=\td_j\,
c^{ij}_m,\quad i,j,m=1,\dots,l+1, \eeq where $c^{ij}_m=\frac{\pal
F^{ij}}{\pal t^m}$. Obviously, the function $F$ satisfies the
equations \beq \frac{\pal^3 F}{\pal t^k \pal t^i \pal
t^j}=\eta_{ij},\quad i,j=1,\dots, l+1 \eeq and the
quasi-homogeneity condition \beq {\mathcal L}_E F=2 F. \eeq From
the properties of a flat pencil of metrics \cite{Du1} it follows
that $F$ also satisfies the associativity equations
\beq\label{zz7}
c_{m}^{ij}\,c_{q}^{m p}=c_{m}^{ip}\,c_{q}^{m j}
\eeq
for any set of fixed indices $i,j,p,q$.
Now the theorem follows from above properties of the function
$F$ and the simple identity $\mathcal{L}_E e=-e$.  The theorem is
proved. \epf


\section{Some examples} \label{sec-5}

In this section we give some examples to illustrate the above
construction of the Frobenius manifold structures. For the sake of simplicity of
notations, instead of $t^1,\dots,t^{l+1}$ we will
redenote the flat coordinates of the metric $\eta^{ij}$ by
$t_1,\dots, t_{l+1}$, and we will also denote $\p_i=\frac{\p}{\p
t_i}$ in the the following examples.

\begin{ex}$[C_3,k=1]$ Let $R$ be the root system of type $C_3$, take $k=1$, then $d_1=d_2=d_3=1$, and
\begin{eqnarray*}
&&y^1={e^{2\,i\pi\,x_{{4}}}} \left( \zeta_{{1}}+\zeta_{{2}}+\zeta_{{3}} \right);\\
&&y^2={e^{2\,i\pi\,x_{{4}}}} \left(
\zeta_{{1}}\zeta_{{2}}+\zeta_{{1}}\zeta_{{3}}+
\zeta_{{2}}\zeta_{{3}}
\right);\\
&&y^3={e^{2\,i\pi\,x_{{4}}}}\zeta_{{1}}\zeta_{{2}}\zeta_{{3}};\\
&&y^4=2\,i\pi\,x_{{4}},
\end{eqnarray*}
where $\zeta_j={e^{2\,i\pi\, \left( x_{{j}}-x_{{j-1}} \right)
}}+{e^{-2\,i\pi\, \left( x_{{j}}-x_{{j-1}} \right) }}+2$ and
$x_0=0$, $j=1,2,3$. The metric $(~,~)^{\sptilde}$ has the form
\begin{equation*}
((dx_i,dx_j)^{\sptilde})=\frac{1}{4\pi^2}\left( \begin
{array}{cccc} ~~1&~~1\,&~~1&~~0\\ ~~1& ~~2&~~2&~~0\\
~~1&~~2&~~3&~~0\\ ~~0&~~0&~~0&-1\end {array}
 \right).
 \end{equation*}
The flat coordinates are
$$t_1=y^{{1}}-2\,e^{y_4},\ t_2=(y^{{2}}-\frac{1}{6}\,y^{{3}})(y^{{3}})^{-\frac14},
\ t_3=(y^{{3}})^{\frac14},\ t_4=y^{{4}}$$
and the intersection form is given by
\begin{eqnarray*}
&&g^{11}=2\,t_{{2}}t_{{3}}\,
{e^{t_{{4}}}}+\frac{1}{3}\,{t_{{3}}}^{4}{e^{t_{{4}}}}+4\,
 {e^{2t_{{4}}}} ;\\
 &&g^{12}=\frac{7}{3}\,{t_{{3}}}^{3}{e^{t_{{4}}}}+\frac{7}{2}\,t_{{2}}{e^{t_{{4}}}};~~
  g^{13}=\frac{5}{2}\,t_{{3}}{e^{t_{{4}}}}; ~~ g^{14}=t_1;\\
&& g^{22}=12\,{t_{{3}}}^{2}{e^{t_{{4}}}}-\frac{1}{4}\,{t_{{2}}}^{2}+\frac{1}{12}\,{t_{{3}}}^{3}t
_{{2}}-{\frac {1}{108}}\,{t_{{3}}}^{6}+\frac{1}{4}\,{\dfrac {{t_{{2}}}^{3}}{{t_
{{3}}}^{3}}};\\
&&g^{23}=2\,t_{{1}}+4\,{e^{t_{{4}}}}-\frac{1}{3}\,t_{{2}}t_{{3}}+{\frac {1}{72}}\,{t_{{
3}}}^{4}-\frac{1}{4}\,{\dfrac {{t_{{2}}}^{2}}{{t_{{3}}}^{2}}};\\
&&g^{24}=\frac{3}{4}\,t_{{2}}; ~~g^{33}=\frac{1}{4}\,{\dfrac
{t_{{2}}}{t_{{3}}}}-\frac{1}{12}\,{t_{{3}}}^{2};~~
g^{34}=\frac{1}{4}\,t_{{3}}; ~~g^{44}=1.
\end{eqnarray*}
The potential has the expression
\begin{eqnarray*}
&&F=\frac{1}{2}\,{t_{{1}}}^{2}t_{{4}}+\frac{1}{2}\,t_{{1}}t_{{2}}t_{{3}}
-{\frac {1}{48}}\,{t_{{ 2}}}^{2}{t_{{3}}}^{2}+{\frac
{1}{1440}}\, t_{{2}}{t_{{3}}}^{5}-{ \frac
{1}{36288}}\,{t_{{3}}}^{8}\\
&&\qquad
+t_{{2}}t_{{3}}{e^{t_{{4}}}}+\frac{1}{6}\,{t_{{3}}}^{4}{e^{t_{{4}}}}+\frac{1}{2}\,
e^{2t_4} +{\frac {1}{48}}\,{\dfrac {{t_{{2}}}^{3}}{t_{{3}}}
}
\end{eqnarray*}
and the Euler vector field is given by
$$E=t_1{\p_1}+\frac{3}{4}t_2{\p_2}+\frac{1}{4}t_3{\p_3}+{\p_4}.$$
\end{ex}

\begin{rem}
If we take $k=2$ or $3$ for the $C_3$ root system, we obtain a Frobenius manifold structure that is isomorphic
to the one given in Example 2.7 [$B_3,k=2$] or  Example 2.8 [$C_3,k=3$] of \cite{DZ1998}.
\end{rem}

\begin{ex}$[C_4,k=1]$
Let $R$ be the root system of type $C_4$. Take $k=1$, then
$d_1=d_2=d_3=d_4=1$, and
\begin{eqnarray*}
&&y^1={e^{2\,i\pi\,x_{{5}}}} \left( \zeta_{{1}}+\zeta_{{2}}+\zeta_{{3}}+\zeta_{{4}} \right);\\
&&y^2={e^{2\,i\pi\,x_{{5}}}} \dsum_{1\leq a<b\leq 4}
\zeta_{{a}}\zeta_{{b}};\\
&&y^3={e^{2\,i\pi\,x_{{5}}}} \dsum_{1\leq a<b<c\leq 4}
\zeta_{{a}}\zeta_{{b}}\zeta_{{c}};\\
&&y^4={e^{2\,i\pi\,x_{{5}}}}\zeta_{{1}}\zeta_{{2}}\zeta_{{3}}\zeta_{{4}}; \\
&&y^5=2\,i\pi\,x_{{5}},
\end{eqnarray*}
where $\zeta_j={e^{2\,i\pi\, \left( x_{{j}}-x_{{j-1}} \right)
}}+{e^{-2\,i\pi\, \left( x_{{j}}-x_{{j-1}} \right) }}+2$ and
$x_0=0$, $j=1,2,3,4$. The metric $(~,~)^{\sptilde}$ has the form
\begin{equation*}
((dx_i,dx_j)^{\sptilde})=\frac{1}{4\pi^2}\left( \begin
{array}{rrrrr}
~~1&~~1&~~1&~~1&~~0\\\noalign{\medskip}1&2&2&2&0\\\noalign{\medskip}1&2&3&3&0\\
\noalign{\medskip}1&2&3&4&0\\ \noalign{\medskip}0&0&0&0&-1\end
{array}
 \right).\end{equation*}
To write down the flat coordinates, we first introduce the variables
\eqa
&&w_1=y^1-2e^{y^5},\
w_2=(y^2-\frac{1}{6}\,y^3+\frac{1}{30}\,y^4)(y^4)^{-\frac16},\nn\\
&& w_3=(y^3-\frac{1}{4}\,y^4) (y^4)^{-\frac23},\quad w_4=(y^4)^{\frac16},\ w_5=y^5.\nn
\eeqa
Then we have
\beq
t_1=w_1,\
t_2=w_2-\frac{1}{12}w_3^2\,w_4,\ t_3=w_3 w_4,\ t_4=w_4,\ t_5=w_5.\nn
\eeq
We omit the presentation of the long expression of the intersection form
and only
write down the potential $F$ here
\begin{eqnarray*}
&&F=\frac
{1}{2}\,{t_{{1}}}^{2}t_{{5}}+{\frac {1}{2}}\,t_{{1}}
t_{{2}}t_{{4}}-{\frac {1}{6912}}
\,{t_{{3}}}^{4}+{\frac{1}{17280}}\,{t_{{3}}}^{3}{t_{{4}}}^{3}\\
&&\qquad -{\frac
{1}{288}}\,t_{{2}}t_{{4}}{t_{{3}}}^{2}-{\frac
{1}{34560}}\,{t_{{4}}}^{ 6}{t_{{3}}}^{2}+{\frac
{1}{24}}\,t_{{1}}{t_{{3}}}^{2}+{\frac
{1}{1440}}\,t_{{3}}{t_{{4}}}^{4 }t_{{2}}\\
&&\qquad -{ \frac
{1}{48}}\,{t_{{2}}}^{2}{t_{{4}}}^{2}-{\frac
{1}{60480}}\,{t_{{4}}}^{7}t_{{2}}+{\frac {1}{
345600}}\,{t_{{4}}}^{9}t_{{3}}-{\frac {1}{
7603200}}\,{t_{{4}}}^{12}\\
&&\qquad+{\frac {1}{12}
}\,{e^{t_{{5}}}}{t_{{3}}}^{2}+\frac
{1}{6}\,{e^{t_{{5}}}}t_{{3}}{t_{{4}}}^{3}+{\frac
{1}{120}}\,{e^{t_{{5}}} }{t_{{4}}}^{6}+
 t_{{2}}t_{{4}}{e^{t_{{5}}}}+\frac {1}{2}\,
{e^{2t_{{5}}}}\\
&&\qquad +{\frac {1}{24}}\,{\dfrac
{t_{{3}}{t_{{2}}}^{2}}{t_{{4}}}}-{\frac {1} {216}}\,{\dfrac
{t_{{2}}{t_{{3}}}^{3}}{{t_{{4}}}^{2}}}+{ \frac
{1}{4320}}\,{\dfrac {{t_{{3}}}^{5}}{{t_{{4}}}^{3}}}.
\end{eqnarray*}
The Euler vector field is given by
$$E=t_1{\p_1}+\frac{5}{6}t_2{\p_2}+\frac{1}{2}t_3{\p_3}+\frac{1}{6}t_4{\p_4}+{\p_5}.$$
\end{ex}

\begin{ex}$[C_4,k=2] \ $Let $R$ be the root system of type $C_4$. Take $k=2$, then
$d_1=1,d_2=d_3=d_4=2$, and
\begin{eqnarray*}
&&y^1={e^{2\,i\pi\,x_{{5}}}} \left( \zeta_{{1}}+\zeta_{{2}}+\zeta_{{3}}+\zeta_{{4}} \right);\\
&&y^2={e^{4\,i\pi\,x_{{5}}}} \dsum_{1\leq a<b\leq 4}
\zeta_{{a}}\zeta_{{b}};\\
&&y^3={e^{4\,i\pi\,x_{{5}}}} \dsum_{1\leq a<b<c\leq 4}
\zeta_{{a}}\zeta_{{b}}\zeta_{{c}};\\
&&y^4={e^{4\,i\pi\,x_{{5}}}}\zeta_{{1}}\zeta_{{2}}\zeta_{{3}}\zeta_{{4}}; \\
&&y^5=2\,i\pi\,x_{{5}},
\end{eqnarray*}
where $\zeta_j$ are defined as in the above example. The metric $(~,~)^{\sptilde}$ has the form
\begin{equation*}
((dx_i,dx_j)^{\sptilde})=\frac{1}{4\pi^2}\left( \begin
{array}{rrrrr}
~~1&~~1&~~1&~~1&~~0\\\noalign{\medskip}1&2&2&2&0\\\noalign{\medskip}1&2&3&3&0\\
\noalign{\medskip}1&2&3&4&0\\
\noalign{\medskip}0&0&0&0&-\frac{1}{2}\end {array}
 \right).\end{equation*}
and the flat coordinates are given by
\eqa
&&t_1=y^{{1}}-4e^{y^5},\ t_2=y^2-2y^1e^{y^5}+6e^{2y^5},\
t_3=(y^3-\frac{1}{6}\,y^4)(y^4)^{-\frac14},\nn\\
&& t_4=(y^4)^{\frac14},\ t_5=y^5.\nn
\eeqa
The intersection form has the following components
\begin{eqnarray*}
&&g^{11}=2t_2-\frac{1}{2}t_1^2+4e^{2t_5};\\
&&g^{12}=6t_1e^{2t_5}+\frac{1}{2}t_4^4e^{t_5}+\frac{1}{2}t_3t_4e^{t_5};\\
&&g^{13}=5\,t_3e^{t_5}+\frac{10}{3}\,t_4^3e^{t_5};
~g^{14}=3\,t_{{4}}{e^{t_{{5}}}};
~g^{15}=\frac{1}{2}t_1;\\
&&g^{22}=2\,{e^{t_{{5}}}}t_{{1}}t_{{3}}t_{{4}}+8\,
{e^{4t_{{5}}}} +8\,
t_{{3}}t_{{4}}{e^{2t_{{5}}}}+\frac{16} {3}\, {e^{2t_{{5}}}}
{t_{{4}}}^{4}+4\, {e^{2t_{ {5}}}}
{t_{{1}}}^{2}+\frac{1}{3}\,{e^{t_{{5}}}}{t_{{4}}}^{4}t_{{1} };
\\
&&g^{23}=\frac{56}{3}\, {e^{2t_{{5}}}} {t_{{4}}}^{3}+7\,   {e^{2t_ {{5}}}}
 t_{{3}}+\frac{7}{3}\,t_{{1}}{t_{{4}}}^{3}{e^{t_{{5}}}}+\frac{7}{2}
\,{e^{t_{{5}}}}t_{{1}}t_{{3}};\\
&&g^{24}=5\, {e^{2t_{{5}}}}t_{{4}}+\frac{5}{2}\,t_{{4}}{e^{t_{{5}}}} t_{{1}} ;~g^{25}=t_2;
\end{eqnarray*}
\begin{eqnarray*}
&&g^{33}=12\,{t_{{4}}}^{2}{e^{t_{{5}}}}t_{{1}}+48\,{t_{{4}}}^{2}
 {e^{ 2t_{{5}}}}
 -\frac{1}{4}\,{t_{{3}}}^{2}+\frac{1}{12}\,t_{{3}}{t_{{4}}}^{3}-\frac{1}{
108}\,{t_{{4}}}^{6}+\frac{1}{4}\,{\dfrac
{{t_{{3}}}^{3}}{{t_{{4}}}^{3}}};\\
&&g^{34}=2\,t_{{2}}+4\,{e^{t_{{5}}}}t_{{1}}+4\, {e^{2t_{{5}}}}
  -\frac{1}{3}\,t_{{4}}t_{{3}}+\frac{1}{72}\,{t_{{4}}}^{4}-\frac{1}{4}\,{\dfrac {{t_{{3}}}^{2}}{{t_{{4}}}^{2}}};\\
&&g^{35}=\frac{3}{4}\,t_{{3}};~g^{44}=\frac{1}{4}\,{\dfrac
{t_{{3}}}{t_{{4}}}}-\frac{1}{12}\,{t_{{4}}}^{2};~g^{45}=\frac{1}{4}t_4;~g^{55}=\frac{1}{2}.
\end{eqnarray*}
The Euler vector field is given by
$$E=\frac{1}{2}t_1{\p_1}+t_2{\p_2}+\frac{3}{4}t_3{\p_3}+\frac{1}{4}t_4{\p_4}+\frac{1}{2}{\p_5}.$$
Finally, we have the potential
\begin{eqnarray*}
&&F=\frac{1}{2}\,{t_{{2}}}^{2}t_{{5}}+\frac{1}{4}\,{t_{{1}}}^{2}t_{{2}}+\frac{1}{2}\,t_{{4}}t_{{3
}}t_{{2}}+{\frac {1}{1440}}\,{t_{{4}}}^{5} t_{{3}}-{\frac
{1}{48}}\,{t_{{4}}}^{2}{t_{{3}}}^{2}\\
&&\qquad-{\frac
{1}{36288}}\,{t_{{4}}}^{8}-{\frac {1}{96}}\,{t_{{1}}}^{4}
+\frac{1}{2}\,{e^{2\,t_{{5}}}}{t_{{1}}}^{2}+\frac{1}{6}\,{e^{t_{{5}}}}t_{{1}}{t_{{4}}}^{
4}+\frac{2}{3}\,{t_{{4}}}^{4}{e^{2\,t_{{5}}}}\\
&&\qquad+{e^{t_{{5}}}}t_{{1}}t_{{3}}t_{{4}}+t_{{3}}t_{{4}}{e^{2\,t
_{{5}}}}+\frac{1}{4}\,{e^{4\,t_{{5}}}}+{\frac {1}{48
}}\,{\dfrac {{t_{{3}}}^{3}}{t_{{4}}}}.
\end{eqnarray*}
\end{ex}
\section{On the Frobenius manifold structures related to the root system of type $B_l$}
\label{sec-6}

For the root system $R$ of type $B_l$, we also choose the standard base $\al_1,\dots, \al_l$
of simple roots as given in \cite{bb}.
As in \cite{DZ1998} we define an indefinite  metric $(~,~)^{\sptilde}$
on $\widetilde V=V\oplus\mathbb R$ such that $\widetilde V$ is the orthogonal direct sum of $V$ and $\mathbb R$,
$V$ endowed with the $W$-invariant Euclidean metric
\beq
(d x_m, dx_n)^{\sptilde}=\frac{1}{4\pi^2}[(1-\frac12\delta_{n,l}) m-\frac{l}4\,\delta_{n,l}\delta_{m,l}],
\quad 1\le m\le n\le l
\eeq
and $\mathbb R$ endowed with the metric
\beq
(dx_{l+1},dx_{l+1})^{\sptilde}=-\frac1{4\pi^2 d_k}.
\eeq
If the label $k$ of the chosen simple root $\al_k$ is less than $l$, then
the numbers $d_j$ are given by
\beq
d_1=1,\dots,d_k=k,\quad d_{k+1}=\dots=d_{l-1}=k, \quad d_l=\frac{k}2.
\eeq
In the case of $k=l$ we have
\beq
d_j=\frac{j}2,\quad j=1,\dots, l-1,\quad d_{l}=\frac{l}4.
\eeq
The basis of the $W_a$-invariant Fourier
polynomials $y_1(\bx),\dots,y_{l-1}(\bx), y_l(\bx)$ are chosen in the same way as we did for the case of $C_l$
\beq\label{zh110}
{y}_j({\bx})=\sigma_j(\zeta_1,\cdots,\zeta_l),\quad j=1,\dots,l-1,
\quad y_l(\bx)=\zeta_1^{\frac12}\dots \zeta_l^{\frac12}
\eeq
with $\zeta_1,\dots,\zeta_{l-1}$ defined by (\ref{zh7b}) and
\beq
\zeta_l=e^{2i \pi (x_{l-1}-2 x_l)}+e^{-2 i\pi (x_{l-1}-2 x_l)}+2.
\eeq
So
\eqa
&&
\zeta_j^{\frac12} = 2 \cos \pi(x_j - x_{j-1}), \quad j\leq l-1
\nn\\
&&
\zeta_l^{\frac12} = 2 \cos \pi(2x_l - x_{l-1}).
\nn
\eeqa
The generators
of the ring $\widetilde{W}^{(k)}(B_l)$ have the same form as
that of (\ref{zh8-b}).
It is easy to see that the components of the resulting metric $(g^{ij}(y))$ coincide
with those corresponding
to the root system of type $C_l$ if we perform the change of coordinates
\beq
y^j\mapsto \bar y^j=y^j,\ y_l(\bx)\mapsto {\bar y}^l=(y^l)^2,\ y^{l+1}\mapsto \bar{y}^{l+1}=y^{l+1},
~~ j=1,\dots,l-1,\nn
\eeq
for $1\le k\le l-1$ and
\beq
y^j\mapsto \bar y^j=y^j,\ y_l(\bx)\mapsto {\bar y}^l=(y^l)^2,\ y^{l+1}\mapsto \bar{y}^{l+1}=\frac12 y^{l+1},
~~j=1,\dots,l-1,\nn
\eeq
for the case when $k=l$. Thus, the Frobenius manifold structure that we obtain in this way from $B_l$,
by fixing the $k$-th vertex of the corresponding Dynkin diagram,
is isomorphic to the one that we obtain from $C_l$ by choosing the $k$-th vertex of the Dynkin diagram of $C_l$.

\section{Concluding remarks}
It remains a challenging problem to
understand whether the constructions of the present paper can be generalized to the root systems of the
types $D_l$, $E_6$, $E_7$, $E_8$, $F_4$, $G_2$ with respect to the choice of an arbitrary vertex on the
Dynkin diagram, as it was suggested in \cite{slodowy} motivating by the results of \cite{wir}. Another
open problem is to obtain an explicit realization of the integrable hierarchies associated with the Frobenius
manifolds of the type $W^{(k)}(R)$. So far this problem was solved only for $k=1$, $R=A_1$, see \cite{CDZ,
DZ2001, DZ04} for details. We plan to study these problems in subsequent publications.



\begin{thebibliography}{99}

\bibitem{CDZ} G. Carlet, B. Dubrovin and Y. Zhang,
The extended Toda hierarchy, Moscow Math. J. {\bf 4} (2004), 313--332.

\bibitem{bb} N. Bourbaki, Groupes et Alg\`ebres de Lie, Chapitres 4, 5
 et 6, Masson, Paris-New York-Barcelone-Milan-Mexico-Rio de Janeiro, 1981.

\bibitem{Du2} B. Dubrovin, Flat pencils of metrics and Frobenius
manifolds, Integrable systems and algebraic geometry (Kobe/Kyoto,
1997), 47--72, World Sci. Publishing, River Edge, NJ, 1998.

\bibitem{Du1}B. Dubrovin, Geometry of 2D topological field
theories, In:Springer Lecture Notes in Math. {\bf 1620}(1996), 120--348.

\bibitem{DZ1998} B. Dubrovin, Y.Zhang, Extended affine Weyl
groups and Frobenius mainfolds, Compositio Mathematica {\bf 111}(1998)
167--219.

\bibitem{DZ2001} B. Dubrovin and Y. Zhang, Normal forms of hierarchies of integrable PDEs,
 Frobenius manifolds and Gromov-Witten invariants, SISSA Preprint 65/2001/FM, arxiv:math.DG/0108160.

\bibitem{DZ04}B. Dubrovin and Y. Zhang,
Virasoro Symmetries of the Extended Toda Hierarchy, Commun. Math. Phys. {\bf 250}(2004), 161--193.

\bibitem{S}K. Saito, On a linear structure of a quotient variety by a finite reflection group,
Publ.RIMS, Kyoto Univ. {\bf 29} (1993), 535--579.

\bibitem{SYS} K. Saito, T. Yano and J. Sekiguchi, On a certain generator
system of the ring of invariants of a finite reflection group,
Comm. Algebra {\bf 8}(4) (1980) 373--408.

\bibitem{slodowy} P. Slodowy, A remark on a recent paper by B. Dubrovin and Y. Zhang, Preprint 1997.

\bibitem{wir} K. Wirthm\"uller,
Torus embeddings and deformations of simple space curves, Acta
Mathematica {\bf 157} (1986), 159-241.



\end{thebibliography}
\end{document}